\documentclass{amsart}
\usepackage{amssymb}
\usepackage{amsmath}
\usepackage{amscd}
\usepackage[matrix,arrow,curve]{xy}

\theoremstyle{definition}

\theoremstyle{remark}

\numberwithin{equation}{section}

\begin{document}

%\title{ Inductive limits of Toeplitz algebras}
%\title{ Inductive sequences and $C^*$-algebras generated by isometries}
%\title{ On inductive sequences and semigroup $C^*$-algebras}
\title[APPROXIMATION OF TENSORS AND  THE TOPOLOGICAL GROUP  STRUCTURE OF   MATRICES]{ A LOW-RANK APPROXIMATION OF TENSORS AND   \\ THE  TOPOLOGICAL  GROUP  STRUCTURE OF INVERTIBLE  MATRICES}
%\rightheadtext{limit automorphisms of Toeplitz algebras}
%    Information for first author
\author{R.N.~Gumerov }
%    Address of record for the research reported here
\address{Chair of Mathematical Analysis, N.I. Lobachevskii Institute of Mathematics and Mechanics,
 Kazan (Volga region) Federal University, Kremlevskaya~35, Kazan,420008, Tatarstan,
 Russian Federation}
%    Current address
%\curraddr{Department of Mathematics and Statistics,
%Case Western Reserve University, Cleveland, Ohio 43403}
\email{Renat.Gumerov@kpfu.ru; rn.gumeroff@gmail.com}
%    \thanks will become a 1st page footnote.
%\thanks{The first author was supported in part by NSF Grant \#000000.}

%    Information for second author
\author{A.S.~Sharafutdinov}
\address{N.I. Lobachevskii Institute of Mathematics and Mechanics,
 Kazan (Volga region) Federal University, Kremlevskaya~35, Kazan,420008, Tatarstan,
 Russian Federation}
%\email{two@maths.univ.edu.au}
%\thanks{Support information for the second author.}

%    General info

%\date{January 1, 1994 and, in revised form, June 22, 1994.}

%\dedicatory{This paper is dedicated to our authors.}


\begin{abstract}
By a tensor we mean an element of a tensor product of  vector spaces over a field. %18
Up to a choice of bases in factors of tensor products,  every tensor may be coordinatized, that is,
represented as an array  consisting of numbers. 
This note is concerned with properties of the~tensor rank that is a natural generalization of the matrix rank. The topological group structure of
invertible matrices is involved in this study. The multilinear matrix multiplication is discussed from a viewpoint of transformation groups.
 We treat a low-rank tensor approximation in finite-dimensional  tensor products. 
 It is shown that the problem on determining a best rank-$n$ approximation for
 a tensor of size $n\times n \times 2$ has no a solution.
 To this end, we make use of an~approximation by matrices with simple spectra.
\end{abstract}

\maketitle

\section*{\textbf{Introduction}}
%As is well known,
%there are close connections between
%topological notions and properties of Banach algebras.
%the  are closely connected with algebraic and topological categories.
%properties of objects
%and  morphisms in the categories of Banach algebras
%the first type categories
% have corresponding analogs for  objects and  morphisms
%in algebraic and topological categories,
%the second type categories,
% and vice versa (see,e.g., \cite[Ch.II]{B}, \cite[Ch.IV]{H}, \cite[Ch.4,7]{P}, \cite[Ch.II-IV]{S71}).
 Tensors are ubiquitous in sciences. The subject of tensors is an active research area in mathematics and its applications
(see, for example, \cite{GumSharHelQ, GumSharHack, GumSharLan} and references therein).

This note is devoted to the tensor rank and a low-rank approximation of tensors.
The~tensor rank can be considered as a measure of complexity of tensors. Therefore,
one is often required to find an approximation of a given tensor by tensors with lower tensor
ranks. In particular, the best low-rank approximation problem for tensors is of great interest in the statistical analysis of multiway data (see, for example,
references in \cite[p.~1085]{GumShardSL}). As is known, in general, the best low-rank approximation problem for tensors is ill-posed \cite{GumSharBini, GumShardSL}.

A part of motivation for this work comes from our study the complexity of tensors  in homological complexes of Banach spaces \cite{GumSharGum1, GumSharGum2}. The main part of motivation comes
from results in \cite{GumShardSL, GumSharTyr, GumSharGumVid} on tensors in finite-dimensional spaces. In this note we consider tensors in finite-dimensional spaces with the Euclidean topology.
Properties of those tensors are closely related to the topological group structure of invertible matrices. Here, we deal with the natural topological group action on a space of tensors. We show
the ill-posedness of the best rank-$n$ approximation problem in the space of tensors of size
$n\times n\times 2$.

The note contains Introduction and two sections. Section~1 contains preliminaries and properties of the tensor rank. We show that the triple consisting of the Cartesian product of general linear groups, the space of tensors and the multilinear matrix multiplication is a topological transformation group. In Section~2 we prove that every element in the tensor space
$\mathbb{C}^{n\times n\times 2}$ may be approximated  arbitrarily closely by elements
whose tensor ranks are equal to $n$.

Introduction and Section~2 were written by R.~N.~Gumerov. Section~1 was written by A.~S.~Sharafutdinov. The results of Section~2 were proved by R.~N.~Gumerov.

%the papers on $C^*$-algebras generated by isometries \cite{C},\cite{C69},\cite{D},\cite{M87}.

%%%%%%%%%%%%%%%%%%%%%%%%%%%%%%%%%%%%%%%%%%%%%%%%%%%%%%%%%%%%%%%%%%%%%%%%

%%%%%%%%%%%%%%%%%%%%%%%%%%%%%%%%%%%%%%%%%%%%%%%%%%%%%%%%%%%%%%%%%%%%%%%%

\section{Tensor rank and its properties}

As usual, $\mathbb{N}$ stands for the set of all natural numbers. In the sequel, $l,m,n\in\mathbb{N}$ and $l,m,n\geqslant 2.$

Throughout the note, $\mathbb{F}$ will denote either the field of complex numbers $\mathbb{C}$
or the field of real numbers $\mathbb{R}$.
For an element $\mathbf{x}\in\mathbb{F}^l$ we use the notation $\mathbf{x}=(x_1,\dots, x_l)^T$, where $x_i\in\mathbb{F}, i=1,\dots, l$.

We denote by $\mathbb{F}^{l\times m}$, or by $M_{l,m}(\mathbb{F})$, the linear space of all matrices $A=(a_{ij})$  of size $l\times m$, where $a_{ij}\in \mathbb{F}$, $i=1,\dots,l, \,j=1,\dots,m.$ The space of all square matrices of order $n$ over the field $\mathbb{F}$ is denoted by $M_{n}(\mathbb{F}).$ The general linear group of degree $n$, that is, the group of invertible matrices in $M_{n}(\mathbb{F})$, is denoted by  $GL_n(\mathbb{F})$. The symbol $E_n$ stands for the identity matrix in
$M_n(\mathbb{F}).$

For $\mathbf{x}=(x_1, \dots, x_l)^T\in\mathbb{F}^l$ and $\mathbf{y}=(y_1, \dots, y_m)^T\in\mathbb{F}^m$, the matrix $\mathbf{x}\otimes \mathbf{y}\in\mathbb{F}^{l\times m}$ is given by
\[
\mathbf{x}\otimes \mathbf{y} = (x_iy_j), \quad \text{where} \quad i=1,\dots,l, \,j=1,\dots,m.
\]

Let $\mathbb{F}^{l\times m\times n}$ be the linear space of all arrays $A=(a_{ijk})$ of size $l\times m\times n$, where $a_{ijk}\in \mathbb{F}$, $i=1,\dots,l, \,j=1,\dots,m, \, k=1,\dots, n.$
For a tensor  $A=(a_{ijk})\in\mathbb{F}^{l\times m\times n}$ we also use the following notation:  \[
A=[A_1|\ldots |A_n],
\]
where, for every $r=1,\dots, n,$ \ the slice $A_r$ is defined by
\[
A_r=(a_{ijr})\in\mathbb{F}^{l\times m},  \quad \text{where} \quad i=1,\dots,l, \,j=1,\dots,m.
\]

For $\mathbf{x}=(x_1, \dots, x_l)^T\in\mathbb{F}^l, \mathbf{y}=(y_1, \dots, y_m)^T\in\mathbb{F}^m, \mathbf{z}=(z_1,\dots,z_n)^T\in\mathbb{F}^n,$ we define the array
$\mathbf{x}\otimes \mathbf{y}\otimes \mathbf{z}\in\mathbb{F}^{l\times m\times n} $ by
\[
\mathbf{x}\otimes \mathbf{y}\otimes \mathbf{z}=(x_iy_jz_k),\quad \text{where} \quad   i=1,\dots,l, j=1,\dots,m,
k=1,\dots, n.
\]

Consider the bilinear mapping $\theta$ and the trilinear mapping $\tau$ defined as follows:
\[
\theta\colon\mathbb{F}^l\times\mathbb{F}^m\longrightarrow\mathbb{F}^{l\times m}: (\mathbf{x},\mathbf{y})\mapsto \mathbf{x}\otimes \mathbf{y};
\]
\[\tau\colon\mathbb{F}^l\times\mathbb{F}^m\times \mathbb{F}^n\longrightarrow\mathbb{F}^{l\times m\times n}:(\mathbf{x},\mathbf{y},\mathbf{z})\mapsto \mathbf{x}\otimes \mathbf{y}\otimes \mathbf{z}.
\]

It is well known that the pairs ($\mathbb{F}^{l\times m},\theta$) and  ($\mathbb{F}^{l\times m\times n},\tau$) are the tensor products for the corresponding
linear spaces. In what follows, elements of the spaces $\mathbb{F}^{l\times m}$ and $\mathbb{F}^{l\times m\times n}$ are called \textit{tensors}.
% is the  tensor  product for $\mathbb{F}^l$ and $\mathbb{F}^m$,
%  and  ($\mathbb{F}^{l\times m\times n},\tau$) is the tensor  product for the  $\mathbb{F}^l, \mathbb{F}^m$ and $ \mathbb{F}^n$.

  For the basics of algebraic tensor products we refer the reader, for instance, to \cite[Part~I, Ch.~3]{GumSharHack},\cite[Ch.~1]{GumSharGreub} and \cite[Ch.~2, \S~7 ]{GumSharHel}.

Both of the $l_1$-norms on the linear spaces $\mathbb{F}^{m\times n}$ and $\mathbb{F}^{l\times m\times n}$ will be denoted by the same symbol $\|\cdot\|_1$. We recall that the value of the~ $l_1$-norm at a tensor $A$ is defined as the sum of absolute values of all entries in $A$.

  All norms on a space of tensors are equivalent and generate the same topology that is called the Euclidean topology. The convergence of a tensor sequence
  \[
  \{A_t\}=\{(a^t_{ijk})\} \subset \mathbb{F}^{l\times m\times n}, \quad t\in \mathbb{N},
   \]
   to a tensor $A=(a_{ijk})\in\mathbb{F}^{l\times m\times n}$ with respect to this topology is exactly the entrywise convergence, that is, for every fixed triple of indices $i=1,\dots,l, j=1,\dots,m,
k=1,\dots, n$, one has the equality
\[
\lim\limits_{t\to +\infty}a^t_{ijk}=a_{ijk}.
\]
In the sequel, we consider spaces of tensors endowed with the Euclidean topologies. The general linear group
$GL_n(\mathbb{F})$ is a topological group with respect to that topology.

\medskip
\textbf{Definition~1.}
	\emph{Tensors $A\in \mathbb{F}^{l\times m}$ and $B\in \mathbb{F}^{l\times m\times n}$ are said to be  \textit{elementary} (or \textit{decomposable}) if $A=\mathbf{a}\otimes \mathbf{b} $ and $B=\mathbf{x}\otimes \mathbf{y}\otimes \mathbf{z}$ for some vectors $\mathbf{a}, \mathbf{x}\in\mathbb{F}^{l},
\mathbf{b},	\mathbf{y}\in\mathbb{F}^{m}$ and
	$\mathbf{z}\in\mathbb{F}^{n}.$}

\medskip
\textbf{Definition~2.}
		\emph{A tensor $A\in \mathbb{F}^{l\times m}$ or a tensor $B\in \mathbb{F}^{l\times m\times n}$ has \textit{the tensor rank}  $r$ if it can be written as a
		sum of $r$ elementary tensors, but no fewer. We will use the notation $rank (A)$ (or  $rank_{\mathbb{F}} (A)$) for the tensor
		rank of $A$. Therefore, we may write
\[
rank(B)= \min \big\{r\mid B=\sum_{i=1}^{r} \mathbf{x}_i\otimes \mathbf{y}_i\otimes \mathbf{z}_i,
\quad \text{where}\quad \mathbf{x}_i\in\mathbb{F}^l, \mathbf{y}_i\in \mathbb{F}^m, \mathbf{z}_i\in\mathbb{F}^n\big\}.
 \]}
	
%\medskip
As is well known, for $A\in \mathbb{F}^{l\times m}$, the tensor rank  $rank (A)$ is exactly the matrix rank and, for $A\in \mathbb{R}^{l\times m}$, the equality $rank_{\mathbb{R}} (A)=rank_{\mathbb{C}} (A)$ is valid. On the other hand, the tensor rank $rank(A)$, where $A\in \mathbb{F}^{l\times m\times n}$, depends on a field $\mathbb{F}$. It is clear,
 that the inequality $rank_\mathbb{C}(A)\leqslant rank_\mathbb{R}(A)$ holds.
\vskip 0.3cm %\cite{GumSharHack}
\textbf{Example.} Let $A=\left[\begin{array}{cc|cc}
1 & 0 & 0& -1\\
0 & -1 &-1& 0\\
\end{array}\right].$ It can be shown  that $rank_{\mathbb{R}}(A)=3$ and $rank_\mathbb{C}(A)=2$
(see \cite[Example~3.44]{GumSharHack}).
\vskip 0.3cm
Further, we introduce the topological group that is the Cartesian product of general linear groups
\[
GL_{l,m,n}(\mathbb{F}):=GL_{l}(\mathbb{F})\times GL_{m}(\mathbb{F})\times GL_{n}(\mathbb{F})
\]
and consider a $GL_{l,m,n}(\mathbb{F})$-action on the space $\mathbb{F}^{l\times m\times n}$
(see also \cite[Section~2.1]{GumShardSL}).
For the notions and facts in the theory of topological transformation groups we refer the reader,
for example, to
\cite{GumSharBredon} and \cite{GumSharVVV}.

Let us take elements \ $A\in\mathbb{F}^{l\times m\times n}$ \ and  \  $ (L,M,N)\in GL_{l,m,n}(\mathbb{F})$ \
 given as follows:
\[
A=(a_{ijk}), \
L=(\lambda_{pi}),\  M=(\mu_{qj}), \ N=(\nu_{rk}).
\]
The tensor $A$ is transformed into the tensor $B = (L,M,N)\cdot A\in\mathbb{F}^{l\times m\times n}$ by the rule:
\[
B=(b_{pqr})\in\mathbb{F}^{l\times m\times n}, \quad \text{where} \quad	
b_{pqr}=\sum_{i,j,k=1}^{l,m,n}\lambda_{pi}\mu_{qj}\nu_{rk}a_{ijk}.
\]

Thus, we have the mapping called \emph{the multilinear matrix multiplication}
\[
\Phi\colon GL_{l,m,n}(\mathbb{F})\times \mathbb{F}^{l\times m\times n}\longrightarrow \mathbb{F}^{l\times m\times n}:((L,M,N),A)\longmapsto(L,M,N)\cdot A,
\]
which was studied in \cite[Sections~2.1, 2.2 and 2.5]{GumShardSL}).

Below the mapping $\Phi$ is considered from a viewpoint of transformation groups.

\vspace{3pt}
\textbf{Proposition~1.} \emph{The following properties are fulfilled:}
\begin{itemize}
\item[$1)$]\, \emph{the triple $\left<GL_{l,m,n}(\mathbb{F}), \mathbb{F}^{l\times m\times n},
\Phi\right>$ is a topological transformation group;}
\item[$2)$]\, \emph{every orbit for the $GL_{l,m,n}(\mathbb{F})$-action consists of elements of the same tensor rank;}
\item[$3)$]\, \emph{the group action of $GL_{l,m,n}(\mathbb{F})$ on $\mathbb{F}^{l\times m\times n}$ is non-effective;}
\item[$4)$]\, \emph{the space $\mathbb{F}^{l\times m\times n}$ is non-homogeneous under the  $GL_{l,m,n}(\mathbb{F})$-action.}
\end{itemize}
\textbf{Proof.}
1)\, We show only the continuity of the multilinear matrix multiplication $\Phi$.

 To this end, we take sequences $\{(L_t, M_t, N_t)\}\subset GL_{l,m,n}(\mathbb{F})$ and $\{A_t\}\subset\mathbb{F}^{l\times m\times n},$ $ t\in \mathbb{N}$,  that converge to $(L,M,N)\in GL_{l,m,n}(\mathbb{F})$ and $A\in\mathbb{F}^{l\times m\times n}$, respectively. Hence, we have the coordinatewise convergence:
 \[
 \lim\limits_{t\to +\infty}L_t=L,\quad \lim\limits_{t\to +\infty}M_t=M,\quad \lim\limits_{t\to +\infty}N_t=N.
 \]

We introduce the following notations:
\[
L_t=(\lambda_{pi}^t),\  M_t=(\mu_{qj}^t), \ N_t=(\nu_{rk}^t), \ A_t=(a_{ijk}^t);
\]
\[
L=(\lambda_{pi}),\  M=(\mu_{qj}), \ N=(\nu_{rk}), \ A=(a_{ijk}).
\]
Then we have the equalities for entries of tensors:
\[
\lim\limits_{t\to\infty}\lambda_{pi}^t=\lambda_{pi}, \
\lim\limits_{t\to\infty}\mu_{qj}^t=\mu_{qj}, \
 \lim\limits_{t\to\infty}\nu_{rk}^t=\nu_{rk}, \
 \lim\limits_{t\to\infty}a_{ijk}^t=a_{ijk}.
\]

Further, we define the constants
 \[
 M_1=\sup\limits_{p,i,t}\left|\lambda_{pi}^t\right|, \ M_2=\sup\limits_{q,j,t}\left|\mu_{qj}^t\right|, \
  M_3=\sup\limits_{r,k,t}\left|\nu_{rk}^t\right|, \ M_4=\sup\limits_{i,j,k,t}\left|a_{ijk}^t\right|.
 \]

In addition, let us set $(L_t,M_t,N_t)\cdot A_t=(b_{ijk}^t)$ and $(L,M,N)\cdot A=(b_{ijk})$.

Then, for every $\varepsilon>0$, there exists $T\in\mathbb{N}$ such that for all $t>T$ we have the~ following inequalities:
\begin{multline*}
\left|b_{pqr}-b_{pqr}^t\right|\leq \left| \sum\limits_{i,j,k=1}^{l,m,n}\lambda_{pi}\mu_{qj}\nu_{rk}a_{ijk}-\sum\limits_{i,j,k=1}^{l,m,n}\lambda_{pi}^t\mu_{qj}^t\nu_{rk}^ta_{ijk}^t\right| \leq\\
\leq\sum\limits_{i,j,k=1}^{l,m,n}\left| \lambda_{pi}\mu_{qj}\nu_{rk}a_{ijk}-\lambda_{pi}^t\mu_{qj}^t\nu_{rk}^ta_{ijk}^t\right|\leq\\
\leq\sum\limits_{i,j,k=1}^{l,m,n}\left(
\left|\lambda_{pi}-\lambda_{pi}^t\right| \left| \mu_{qj} \nu_{rk} a_{ijk}\right|+
\left| \mu_{qj}-\mu_{qj}^t\right| \left|\lambda_{pi}^t \nu_{rk} a_{ijk}\right|\right)+\\
+
\sum\limits_{i,j,k=1}^{l,m,n}\left(\left| \nu_{rk}-\nu_{rk}^t\right|\left|\lambda_{pi}^t \mu_{qj}^t a_{ijk}\right|+
\left| a_{ijk}-a_{ijk}^t\right|\left|\lambda_{pi}^t \mu_{qj}^t \nu_{rk}^t\right|
\right) \leq\\
\leq \sum\limits_{i,j,k=1}^{l,m,n} \left(\frac{\varepsilon M_2M_3M_4}{4lmnM_2M_3M_4}+\frac{\varepsilon M_1M_3M_4}{4lmnM_1M_3M_4}+\frac{\varepsilon M_1M_2M_4}{4lmnM_1M_2M_4}+\frac{\varepsilon M_1M_2M_3}{4lmnM_1M_2M_3}\right)\leq\\
\leq \sum\limits_{i,j,k=1}^{l,m,n}\frac{\varepsilon}{lmn}\leq\varepsilon .
\end{multline*}
Hence, the sequence $\Phi((L_t,M_t,N_t), A_t)$ converges to the tensor $\Phi((L,M,N), A)$, as required.

2)\, See the proof of Lemma~2.3(2) in \cite{GumShardSL}.

3)\,Take $(L,M,N)=(\alpha E_l, \beta E_m, \gamma E_n)$
with arbitrary scalars $\alpha, \beta, \gamma\in\mathbb{F}$ satisfying the~condition $\alpha\beta\gamma=1.$
Then, for every $A=(a_{ijk})\in\mathbb{F}^{l\times m\times n},$ we have
\[
\Phi((L,M,N),A)=(\alpha\beta\gamma a_{ijk})=(a_{ijk})=\Phi((E_l,E_m,E_n),A).
\]
This shows that the action is non-effective.

4)\,  Consider tensors $A=\mathbf{x}_1\otimes\mathbf{y}_1\otimes\mathbf{z}_1$ and $B=\mathbf{x}_1\otimes\mathbf{y}_1\otimes\mathbf{z}_1+\mathbf{x}_2\otimes\mathbf{y}_2\otimes\mathbf{z}_2,$ where $\{\mathbf{x}_1,\mathbf{x}_2\}\subset \mathbb{F}^l, \{\mathbf{y}_1,\mathbf{y}_2\}\subset \mathbb{F}^m$ and $ \{\mathbf{z}_1,\mathbf{z}_2\}\subset \mathbb{F}^n$ are pairs of linear independent vectors. Then, one has $rank\,(A)=1$ and $rank\,(B)=2$. Using item~2), we obtain the desired conclusion. $\Box$

\vspace{5pt}
We have the following proposition on the semicontinuity of the tensor rank (see \cite[Proposition~4.3, Theorem~4.10]{GumShardSL}).

\bigskip
\textbf{Proposition~2.}\emph{ The following properties are fulfilled:
\begin{itemize}
\item[$1)$]\, for every $r\leq\min(l,m)$, the set $\mathcal{S}_r(l,m):=\{A\in\mathbb{F}^{l\times m}| rank(A)\leq r\}$ is closed;
\item[$2)$]\, there exists $r$ such that the set $\mathcal{S}_r(l,m,n):=\{B\in\mathbb{F}^{l\times m\times n}| rank(B)\leq r\}$ is not closed.%или лучше for all r\geq2?
\end{itemize}}

%\vskip 0.3cm

\section{The topological group of invertible matrices and approximations of matrices and  tensors}

In this section we consider a low-rank approximation of tensors in the space $\mathbb{C}^{n\times n\times 2}$.

To this end, first, we  formulate Bi's criterion for square-type tensors in the space $\mathbb{C}^{m\times m\times n}$ (see \cite[Proposition~2.5] {GumSharSSM}).

\vspace{3pt}
\textbf{Proposition~3.}
		\emph{Let $A=[A_1|\ldots |A_n]$ be a tensor in $\mathbb{C}^{m\times m\times n}$, where $n\geqslant 2$, and
let  $A_1\in \mathbb{C}^{m\times m} $ be a nonsingular matrix. Then %$rank_{\otimes}A=n$
%$trank(A)=n$
the tensor rank of $A$ is equal to $m$ if and only if the matrices
		$A_2A_1^{-1},\ldots, A_nA_1^{-1}$ can be diagonalized simultaneously.}

Second, we recall some algebraic and topological definitions and facts about square matrices.

\emph{An eigenvalue} of a matrix is said to be \emph{simple} if its algebraic
multiplicity equals one. \emph{The spectrum} of a matrix is said to
be \emph{simple} provided that all eigenvalues of a given matrix are
simple. In other words, if all eigenvalues are pairwise distinct.
Certainly,  a matrix with a simple spectrum is diagonalizable.

We recall that a mapping $f:X \rightarrow Y$ between two topological
spaces is said to be \textit{open}, if for any open set $O$ in $X$
the image $f(O)$ is open in $Y$. For example, if a mapping
$f:GL_n(\mathbb{C})\rightarrow GL_n(\mathbb{C})$ is
a surjective continuous homomorphism then it is
open \cite[Theorem~5.29]{GumSharHR}.

Making use of the topological group structure of the general linear group $GL_n(\mathbb{C})$, one can prove the following statement
\cite[Proposition~4]{GumSharGumVid}.

\vspace{3pt}
\textbf{Proposition~4.}\emph{
Let $f_1,f_2,\ldots, f_k:GL_n(\mathbb{C})\longrightarrow GL_n(\mathbb{C})$   be
a finite family of self-mappings of the general linear group.
Assume that at least one of these
mappings is open with respect to the~Euclidean topology. Let $A_1,A_2,\ldots, A_k$ be  arbitrary
matrices in $M_n(\mathbb{C})$ and let $\|\cdot\|$ be a~norm on $M_n(\mathbb{C})$. Then for every $\varepsilon > 0$ there
exists a finite family  $A_{1\varepsilon},A_{2\varepsilon},\ldots,
A_{k\varepsilon}$ consisting of  invertible matrices with simple
spectra such that the inequalities
\[
\|A_1-A_{1\varepsilon}\|<\varepsilon, \,
\|A_2-A_{2\varepsilon}\|<\varepsilon,
\ldots,\|A_k-A_{k\varepsilon}\|<\varepsilon
\]
hold and the product matrix
\[
f_1(A_{1\varepsilon})f_2(A_{2\varepsilon})\ldots
f_k(A_{k\varepsilon})
\]
 has a simple spectrum.}

\vspace{3pt}
For the case of two mappings, we put $f_1$ and $f_2$ to be  the identity mapping
and the inverse mapping respectively:
\[
f_1: GL_n(\mathbb{C})\longrightarrow GL_n(\mathbb{C}): X\longmapsto X;
\qquad
f_2: GL_n(\mathbb{C})\longrightarrow GL_n(\mathbb{C}): X\longmapsto X^{-1}.
\]
Obviously, both  of these mappings are open  with respect to the Euclidean topology in $GL_n(\mathbb{C})$. Therefore, as a consequence of the preceding proposition, we have

\vspace{3pt}
\textbf{Corollary~1.}
\textit{Let $A$ and $B$ be
matrices $M_n(\mathbb{C})$ and let $\|\cdot\|$ be a norm on $M_n(\mathbb{C})$. Then for every $\varepsilon > 0$ there
exists a pair of  matrices $A_{\varepsilon}$ and
$B_{\varepsilon}$ in  $GL_n(\mathbb{C})$ with simple spectra such that the  inequalities
\[
\|A-A_{\varepsilon}\|<\varepsilon \quad \text{and} \quad \|B-B_{\varepsilon}\|<\varepsilon
\]
hold and the product matrix $A_{\varepsilon}B_{\varepsilon}^{-1}$ has a
simple spectrum.}
\vskip 0.3cm

It is worth noting that one can use Corollary~1 for estimating
the tensor rank of inverse matrices in the case when given matrices  are the factors of the Kronecker  products (see \cite{GumSharTyr}).

We make use of the above-mentioned results to prove the following assertion.

\bigskip
\textbf{Proposition~5.}~\label{propos5}
\emph{Let $A$ be a tensor in $ \mathbb{C}^{n\times n\times 2}$ and let $\|\cdot\|$ be a norm on $ \mathbb{C}^{n\times n\times 2}$. Then the
equality
\[
		\inf \big\{\|A- B\|: B\in \mathbb{C}^{n\times n\times 2}  \quad \text{and} \quad  \text{rank}\,(B)=n\big\}=0
		\]
holds, that is,
the tensor $A$ may be approximated by tensors whose tensor ranks are equal to $n$.}

\textbf{Proof.}
We set $A=[A_1|A_2],$  where $A_1$ and $A_2$ are square matrices of size $n\times n$.

Let us fix $\varepsilon >0 $. Using Corollary~1, we take two invertible matrices $B_{1}$ and $B_{2}$ of size $n\times n$
such that the  inequalities%satisfying the following conditions:
\[
\|A_1-B_1\|_1< \frac{\varepsilon}{2} \quad  \text{and}\quad  \|A_2-B_2\|_1< \frac{\varepsilon}{2}
\]
hold and the product matrix $B_2B_1^{-1}$ has a simple spectrum.

Consider the tensor $B=[B_1|B_2]$.  By Bi's criterion, since the matrix $B_2B_1^{-1}$ is diagonalizable, the tensor rank of $B$ is equal to $n$. Moreover,  we have
the following estimation:
\[
\|A-B\|_1=\|A_1-B_1\|_1+\|A_2-B_2\|_1 < \varepsilon.
\]

In view of the equivalence of the norms $\|\cdot\|$ and $\|\cdot\|_1$, the rest is clear. $\Box$

\vspace{3pt}
It is known (see \cite{GumerovJaJa}, \cite[Theorem~4.3]{GumerovSharSSMlaa}
), that the maximal value of the tensor rank
on the space $\mathbb{C}^{n\times n\times 2}$ is given by
\[
mrank\,(n,n,2):=\max \left\{ rank\,(A) \, | \, A\in \mathbb{C}^{n\times n\times 2} \right \}=
n+ \left\lfloor \frac{n}{2}\right\rfloor,
\]
where the symbol $\lfloor \cdot \rfloor$ means the integer part of a real number.
Therefore $mrank\,(2n,2n,2)=3n$ for every $n\in \mathbb{N}.$
This fact together with the Proposition~5 guarantees that, generally speaking, the tensor rank can leap an arbitrary large gap (see also \cite[Section~4.5]{GumShardSL}). More precisely, we have

\vspace{3pt}
\textbf{Corollary~2.}
\textit{Let $n\in \mathbb{N}.$ There exists a tensor $A\in \mathbb{C}^{2n\times 2n\times 2}$ with $rank\,(A) =3n$ and a sequence of tensors $\{A_k\} \subset \mathbb{C}^{2n\times 2n\times 2}$, $k\in \mathbb{N},$ such that $rank\,(A_k)=2n $ for every $k\in \mathbb{N}$ and
\[
\lim_{k\rightarrow +\infty} A_k =A,
\]
where the limit is taken in the Euclidean topology.
}

Finally,  we can conclude that the tensor rank is not  semicontinuous on the tensor space $\mathbb{C}^{n\times n\times 2}$ endowed with the Euclidean topology.

\vspace{3pt}
\textbf{Corollary~3.}
\textit{Let $n\geqslant 2$. In the tensor space $\mathbb{C}^{n\times n\times 2}$ endowed with the Euclidean topology the set of tensors
		\[
		\big\{\,T\in \mathbb{C}^{n\times n\times 2} \, | \, rank \,(T)\leqslant n\,\big\}
		\]
is not closed and its closure coincides with the whole space $\mathbb{C}^{n\times n\times 2}$.}
%\vskip 0.3cm

\section*{Acknowledgments}

The authors are grateful to participants of the seminar on functional and numerical analysis \
"Tensor Analysis" \ at Kazan Federal University, the seminar \ "Quantum Functional Analysis and Its Applications" \  at Kazan State Power Engineering University and International Conference \ "Probability Theory and Mathematical Statistics" (November 7--10, 2017, Kazan) for helpful discussions of the~content of this note.

%\centerline{Acknowledgment}
%The author thanks Professor S.A.Grigorian
%who called the author's attention to the papers on Toeplitz algebras and
%for the helpful discussions of the results of this article.
\bibliographystyle{amsplain}

\end{document}